\documentclass[10 pt, a4 paper, twoside]{article}
 \usepackage{amssymb,amsmath,graphicx,enumerate}
 \usepackage{longtable,mathrsfs}
 \usepackage{hyperref}
  \numberwithin{equation}{section}
 \newtheorem{thm}{Theorem}[section]
 
 \newtheorem{lem}{Lemma}[section]

 \newtheorem{defn}{Definition}[section]

 \allowdisplaybreaks
 \begin{document}
\begin{center}\textbf{HARTMAN-WINTNER-TYPE INEQUALITY FOR FRACTIONAL DIFFERENTIAL EQUATIONS WITH PRABHAKAR DERIVATIVE}\end{center}
\begin{center}
\textbf{ Deepak B. Pachpatte${^{1}}$, Narayan G. Abuj${^{2}}$ and Amol D. Khandagale${^{3}}$}\\
Department of mathematics,\\ Dr. Babasaheb Ambedkar Marathwada University,\\ Aurangabad-431004 (M.S.) India.\\
{ ${{^1}}$pachpatte@gmail.com,${{^2}}$abujng@gmail.com and ${{^3}}$kamoldsk@gmail.com}
\end{center}
\begin{abstract}
 In this paper, we consider a nonlocal fractional boundary value problem with Prabhakar derivative and obtained a Hartman-Wintner type inequality for it.
\end{abstract}
${\textbf{\emph{Mathematics Subject Classification}}}: 26A33; 33E12; 34A08; 34A40$. \\
\textbf{Keywords:}  \emph{Hartman-Wintner-type inequality; Fractional boundary value problem, nonlocal boundary conditions; Prabhakar derivative}.
 \section{Introduction}
 \paragraph{}  Hartman and Wintner  \cite{HW} considered the following boundary value problem
 \begin{equation}\label{1}
 \left\{ {\begin{array}{*{20}{l}}{x^{\prime\prime}(t)}+q(t)x(t)= 0,\quad a<t< b,\quad\\
 x(a)=x(b)=0, \end{array}} \right.
 \end{equation}
 and proved the inequality if \eqref{1}have a nontrivial solution
 \begin{equation}\label{222}
 \int_{a}^{b}(b-s)(s-a)q^{+}(s)ds >{b-a},
 \end{equation}
 provided \eqref{1} has a nontrivial solution, where $q^{+}(s)= max\{q(s),0\}$.\\
 \paragraph{} Lyapunov  \cite{AML}  proved that if \eqref{1} have a nontrivial solution, then
 \begin{equation}\label{111}
   \int_{a}^{b}|q(s)|ds>\frac{4}{b-a}.
 \end{equation}
 This inequality \eqref{111} is useful in various branches of mathematics such as oscillation theory, disconjugacy and eigenvalue problems.
 The Lyapunov inequality \eqref{111} can be deduced from \eqref{222}  using the fact that
 \begin{equation}
\mathop{\max}\limits_{s\in[a,b]}(b-s)(s-a)= \frac{(b-a)^{2}}{4}.\\
\end{equation}
\paragraph{} The generalizations and extensions of the Lyapunov inequality \eqref{111} exist in the literature \cite{CS1,BH, CD,CS2,PB,LYH,PP,TUC}. Recently, some Lyapunov type inequalities
were obtained for different fractional boundary value problem using various differential operators \cite{FR1,RAC2, MJ1, MJ2, RB, RS, DNA,ND}.
\paragraph{} Cabrera and et al. \cite{CSS} have considered the following nonlocal fractional boundary value problem
 \begin{equation}\label{3}
 \left\{ {\begin{array}{*{20}{l}}{D}_{a}^{\alpha}x(t)+q(t)x(t)= 0,\quad a<t<b,\\
 x(a)=x^{\prime}(a)=0, x^{\prime}(b)=\beta x(\xi), \end{array}}\right.
 \end{equation}
 where $D_{a}^{\alpha}$ denotes the standard Riemann-Liouville fractional derivative of order $\alpha$, $a<\xi<b$, $0\leq\beta(\xi-a)^{\alpha-1}<(\alpha-1)(b-a)^{\alpha-2}$, $q(t)$ is
 continuous real valued function on $[a,b]$. The authors obtained the Hartman-Wintner-type inequality
 \begin{equation}\label{2}
 \int_{a}^{b}(b-s)^{\alpha-2}(s-a)|q(s)|ds \geq\bigg(1+\frac{\beta(b-a)^{\alpha-1}}{(\alpha-1)(b-a)^{\alpha-2}-\beta(\xi-a)^{\alpha-1}}\bigg)^{-1}\Gamma(\alpha).
 \end{equation}
\paragraph{} Motivated by above work, in this paper we consider the following nonlocal fractional boundary problem
  \begin{equation}\label{3}
 \left\{ {\begin{array}{*{20}{l}}\textbf{D}_{\rho,\mu,\omega,a^{+}}^{\gamma}x(t)+q(t)x(t)= 0,\quad a<t<b,\quad 2<\mu\leq3,\\
 x(a)=x^{\prime}(a)=0, x^{\prime}(b)=\beta x(\xi), \end{array}}\right.
 \end{equation}
 where $\textbf{D}_{\rho,\mu,\omega,a^{+}}^{\gamma}$ denotes the Prabhakar derivative of order $\mu$.\\
 $a<\xi<b, 0\leq\beta(\xi-a)^{\mu-1}<(\mu-1)(b-a)^{\mu-2}$, $q:[a,b]\rightarrow\mathbb{R}$ is real valued continuous function and obtained the Hartman-Wintner-type inequality for it.
\section{Preliminaries} In this section, we give some basic definitions and lemmas that will be important to us in the sequel.
\begin{defn}
\cite{TRP} The generalized Mittag-Leffler function with three parameters is defined as,
\begin{equation}\label{a1}
E_{\rho,\mu}^{\gamma}(z)= \sum_{k=0}^\infty\frac{(\gamma)_kz^k}{\Gamma(\rho k+\mu)k!},\qquad \gamma, \rho, \mu \in\mathbb{C}, \Re(\rho)>0,
\end{equation}
where ${(\gamma)_k}$ is Pochhammer symbol defined by,\\
${(\gamma)_0}=1,\quad {(\gamma)_k}={\gamma(\gamma+1)...(\gamma+k-1)}$, for  $k=1,2,...$.\\
For $\gamma=1$, the generalized Mittag-Leffler function \eqref{a1} reduces to the two-parameter Mittag-Leffler function given by
\begin{equation}\label{a2}
E_{\rho,\mu}(z):=E_{\rho,\mu}^{1}(z)=\sum_{k=0}^\infty\frac{z^k}{\Gamma(\rho k+\mu)},\qquad\rho, \mu \in \mathbb{C},\quad \Re(\rho)>0,
\end{equation}
and for $\mu=\gamma=1$, this function coincides with the classical Mittag-Leffler function $E_{\rho}(z)$
\begin{equation}\label{a3}
E_{\rho}(z):=E_{\rho, 1}^{1}(z)=\sum_{k=0}^\infty\frac{z^k}{\Gamma(\rho k+1)}, \qquad \rho\in \mathbb{C},\quad \Re(\rho)>0.
\end{equation}
Also, for $\gamma=0$ we have
$E_{\rho,\mu}(z)=\frac{1}{\Gamma(\mu)}$.
\end{defn}
 \begin{defn}
\cite{RG} Let $f\in L^{1}[0,b]$, $0<x<b\leq\infty$, the prabhakar integral operator including generalized Mittag-Leffler function \eqref{a1} is defined as follows
\begin{equation}\label{a4}
\textbf{E}_{\rho,\mu,\omega,0+}^{\gamma}f(x)dx =\int_{0}^{x}(x-u)^{\mu-1}E_{\rho,\mu}^{\gamma}(\omega(x-u)^{\rho})f(u)du,\quad x>0
\end{equation}
where ${\rho,\mu,\omega,\gamma}\in \mathbb{C},$ with $\Re(\rho), \Re(\mu)>0.$\\
If for $\gamma=0$, the prabhakar integral operator coincides with the Riemann-Liouville fractional integral of order $\mu$;
$$\textbf{E}_{\rho,\mu,\omega,0+}^{0}f(x)=I_{0+}^{\mu}f(x),$$
where the Riemann-Lioville fractional integral is defined as
\begin{equation}\label{a5}
I_{0+}^{\mu}f(x)=\frac{1}{\Gamma(\mu)}\int_{0}^{x}(x-t)^{\mu-1}f(t)dt, \quad \mu\in \mathbb{C}, \Re(\mu)>0.
\end{equation}
\end{defn}
\begin{defn}
\cite{RG}  Let $f\in L^{1}[0,b]$, $0<x<b\leq\infty$, the Prabhakar derivative is defined as
\begin{equation}\label{a6}
\textbf{D}_{\rho,\mu,\omega,0+}^{\gamma}f(x)=\frac{d^m}{d{x}^m}\textbf{E}_{\rho, m-\mu,\omega,{0+}}^{-\gamma}f(x),
\end{equation}
where $\rho,\mu,\omega,\gamma\in \mathbb{C},$ with $\Re(\rho)>0$, $\Re(\mu)>o$, $m-1<\Re(\mu)<m$.\\
We note that the Prabhakar derivative genralizes the Riemann-Liouville fractional derivative
\begin{equation}\label{a7}
D_{0+}^{\mu}f(x)=\frac{d^m}{d{x}^m}\bigg(I_{0+}^{m-\mu}f \bigg)(x),\quad\mu \in\mathbb{C}, \Re(\mu)>0, m-1<\Re(\mu)<m.
\end{equation}
\end{defn}
\begin{lem}\label{L1} \cite{EA} If $f(x)\in C(a,b)\bigcap L(a,b),$ then
\begin{equation}\label{b3}
 \textbf{D}_{\rho,\mu,\omega,a+}^{\gamma}\textbf{E}_{\rho,\mu,\omega,a+}^{\gamma}f(x)=f(x),
\end{equation}
and if\\ $f(x),\textbf{D}_{\rho,\mu,\omega,a+}^{\gamma}f(x)\in C(a,b)\bigcap L(a,b)$ \text{then for} $c_{j} \in \mathbb{R},$ and $m-1<\mu\leq m$,\\
we have
\begin{align}\label{bb2}
\textbf{E}_{\rho,\mu,\omega,a+}^{\gamma} \textbf{D}_{\rho,\mu,\omega,a+}^{\gamma}f(x)&=f(x)+c_{1}(x-a)^{\mu-1}E_{\rho,\mu}^{\gamma}(\omega(x-a)^\rho){\nonumber}\\
&\quad +c_{2}(x-a)^{\mu-2}E_{\rho,\mu-1}^{\gamma}(\omega(x-a)^\rho)+...{\nonumber}\\
&\quad +c_{m}(x-a)^{\mu-m}E_{\rho,\mu-m+1}^{\gamma}(\omega(x-a)^\rho).
\end{align}
\end{lem}
\section{Main Results}
 \begin{thm}\label{T1}
 Assume that $2<\mu\leq3$ and $x\in C[a,b]$. If the nonlocal fractional boundary value problem \eqref{3} has unique nontrivial solution, then it satisfies
 \begin{align*}
&x(t)=\int_{a}^{b}G(t,s)x(s)ds
+\frac{\beta(t-a)^{\mu-1}E_{\rho,\mu}^{\gamma}(\omega(t-a)^{\rho})}{(b-a)^{\mu-2}E_{\rho,\mu-1}^{\gamma}(\omega(b-a)^{\rho})-\beta(\xi-a)^{\mu-1}E_{\rho,\mu}^{\gamma}(\omega(\xi-a)^{\rho})}\\
&\qquad\qquad\qquad\qquad\qquad\qquad\times\int_{a}^{b}G(\xi,s)x(s)ds,
\end{align*}
where the Green's function is defined as
\begin{equation}\label{b4}
G(t,s)= \left\{
{\begin{array}{*{20}{l}}\frac{(t-a)^{\mu-1}E_{\rho,\mu}^{\gamma}(\omega(t-a)^\rho){(b-s)^{\mu-2}E_{\rho,\mu-1}^{\gamma}(\omega(b-s)^\rho)}}{{(b-a)^{\mu-2}E_{\rho,\mu-1}^{\gamma}(\omega(b-a)^\rho)}}\\
 -(t-s)^{\mu-1}E_{\rho,\mu}^{\gamma}(\omega(t-s)^\rho),\quad\quad\quad\quad\quad\quad\quad a\leq s\leq t\leq b, \\
 \frac{(t-a)^{\mu-1}E_{\rho,\mu}^{\gamma}(\omega(t-a)^\rho){(b-s)^{\mu-2}E_{\rho,\mu-1}^{\gamma}(\omega(b-s)^\rho)}}{{(b-a)^{\mu-2}E_{\rho,\mu-1}^{\gamma}(\omega(b-a)^\rho)}},\quad a\leq t\leq s\leq b. \end{array}} \right.
 \end{equation}
 \end{thm}
\textbf{Proof}. From lemma \ref{L1}, the general solution to \eqref{3} in $C[a,b]$ can be written as follows
\begin{align}
x(t)= \nonumber &c_{1}(t-a)^{\mu-1}E_{\rho,\mu}^{\gamma}(\omega(t-a)^\rho)+c_{2}(t-a)^{\mu-2}E_{\rho,\mu-1}^{\gamma}(\omega(t-a)^\rho)\\\nonumber
&+c_{3}(t-a)^{\mu-3}E_{\rho,\mu-2}^{\gamma}(\omega(t-a)^\rho)-\int_{a}^{t}(t-s)^{\mu-1}E_{\rho,\mu}^{\gamma}(\omega(t-s)^\rho)q(s)x(s)ds.\\
\end{align}
Employing the first boundary condition $x(a)= x^{\prime}(a)=0$ we obtain $c_{2}= c_{3}=0$.
Therefore
\begin{align}
x(t)= c_{1}(t-a)^{\mu-2}E_{\rho,\mu}^{\gamma}(\omega(t-a)^\rho)-\int_{a}^{t}(t-s)^{\mu-1}E_{\rho,\mu}^{\gamma}(\omega(t-s)^\rho)q(s)x(s)ds. \nonumber\\
\end{align}
For second boundary condition we obtain,
\begin{align*}
x^{\prime}(t)=c_{1}(t-a)^{\mu-2}E_{\rho,\mu-1}^{\gamma}(\omega(t-a)^\rho)
-\int_{a}^{t}(t-s)^{\mu-2}E_{\rho,\mu-1}^{\gamma}(\omega(t-s)^\rho)q(s)x(s)ds.
\end{align*}
Employing the second boundary condition $x^{\prime}(b)= \beta x(\xi)$ we get
\begin{align*}
&c_{1}(b-a)^{\mu-2}E_{\rho,\mu-1}^{\gamma}(\omega(b-a)^\rho)-\int_{a}^{b}(b-s)^{\mu-2}E_{\rho,\mu-1}^{\gamma}(\omega(b-s)^\rho)q(s)x(s)ds=\\
&\beta c_{1}(\xi-a)^{\mu-1}E_{\rho,\mu}^{\gamma}(\omega(\xi-a)^\rho)-\beta \int_{a}^{\xi}(\xi-s)^{\mu-1}E_{\rho,\mu}^{\gamma}(\omega(\xi-s)^\rho)q(s)x(s)ds,\\
&\Rightarrow c_{1}(b-a)^{\mu-2}E_{\rho,\mu-1}^{\gamma}(\omega(b-a)^\rho)-\beta c_{1}(\xi-a)^{\mu-1}E_{\rho,\mu}^{\gamma}(\omega(\xi-a)^\rho)= \\
&\int_{a}^{b}(b-s)^{\mu-2}E_{\rho,\mu-1}^{\gamma}(\omega(b-s)^\rho)q(s)x(s)ds - \beta \int_{a}^{\xi}(\xi-s)^{\mu-1}E_{\rho,\mu}^{\gamma}(\omega(\xi-s)^\rho)q(s)x(s)ds,\\
&\Rightarrow c_{1}\bigg[(b-a)^{\mu-2}E_{\rho,\mu-1}^{\gamma}(\omega(b-a)^\rho)-\beta(\xi-a)^{\mu-1}E_{\rho,\mu}^{\gamma}(\omega(\xi-a)^\rho)\bigg]=\\
&\int_{a}^{b}(b-s)^{\mu-2}E_{\rho,\mu-1}^{\gamma}(\omega(b-s)^\rho)q(s)x(s)ds- \beta \int_{a}^{\xi}(\xi-s)^{\mu-1}E_{\rho,\mu}^{\gamma}(\omega(\xi-s)^\rho)q(s)x(s)ds, \\
&\Rightarrow c_{1}=\frac{1}{(b-a)^{\mu-2}E_{\rho,\mu-1}^{\gamma}(\omega(b-a)^\rho)- \beta (\xi-a)^{\mu-1}E_{\rho,\mu}^{\gamma}(\omega(\xi-a)^\rho)}\\
&\qquad\qquad\qquad\times\int_{a}^{b}(b-s)^{\mu-2}E_{\rho,\mu-1}^{\gamma}(\omega(b-s)^\rho)q(s)x(s)ds\\
&\qquad\qquad-\frac{\beta}{(b-a)^{\mu-2}E_{\rho,\mu-1}^{\gamma}(\omega(b-a)^\rho)- \beta (\xi-a)^{\mu-1}E_{\rho,\mu}^{\gamma}(\omega(\xi-a)^\rho)}\\
&\qquad\qquad\qquad\times\int_{a}^{\xi}(\xi-s)^{\mu-1}E_{\rho,\mu}^{\gamma}(\omega(\xi-s)^\rho)q(s)x(s)ds.
\end{align*}
Thus the solution $x(t)$ becomes
\begin{align*}
x(t)=&\frac{(t-a)^{\mu-1}E_{\rho,\mu}^{\gamma}(\omega(t-a)^\rho)}{(b-a)^{\mu-2}E_{\rho,\mu-1}^{\gamma}(\omega(b-a)^\rho)- \beta (\xi-a)^{\mu-1}E_{\rho,\mu}^{\gamma}(\omega(\xi-a)^\rho)}\\
&\quad\times\int_{a}^{b}(b-s)^{\mu-2}E_{\rho,\mu-1}^{\gamma}(\omega(b-s)^\rho)q(s)x(s)ds \\
&-\frac{\beta (t-a)^{\mu-1}E_{\rho,\mu}^{\gamma}(\omega(t-a)^\rho)}{(b-a)^{\mu-2}E_{\rho,\mu-1}^{\gamma}(\omega(b-a)^\rho)- \beta (\xi-a)^{\mu-1}E_{\rho,\mu}^{\gamma}(\omega(\xi-a)^\rho)}\\
&\quad\times\int_{a}^{\xi}(\xi-s)^{\mu-1}E_{\rho,\mu-1}^{\gamma}(\omega(\xi-s)^\rho)q(s)x(s)ds\\
&-\int_{a}^{t}(t-s)^{\mu-1}E_{\rho,\mu}^{\gamma}(\omega(t-s)^\rho)q(s)x(s)ds.
\end{align*}
Taking into account that
\begin{align*}
&\frac{E_{\rho,\mu}^{\gamma}(\omega(t-a)^\rho)}{(b-a)^{\mu-2}E_{\rho,\mu-1}^{\gamma}(\omega(b-a)^\rho)- \beta (\xi-a)^{\mu-1}E_{\rho,\mu}^{\gamma}(\omega(\xi-a)^\rho)} \\
&=\frac{E_{\rho,\mu}^{\gamma}(\omega(t-a)^\rho)}{{(b-a)^{\mu-2}E_{\rho,\mu-1}^{\gamma}(\omega(b-a)^\rho)}}\bigg(\frac{(b-a)^{\mu-2}E_{\rho,\mu-1}^{\gamma}(\omega(b-a)^\rho)}{(b-a)^{\mu-2}E_{\rho,\mu-1}^{\gamma}(\omega(b-a)^\rho)- \beta (\xi-a)^{\mu-1}E_{\rho,\mu}^{\gamma}(\omega(\xi-a)^\rho)}\bigg)\\
&=\frac{E_{\rho,\mu}^{\gamma}(\omega(t-a)^\rho)}{{(b-a)^{\mu-2}E_{\rho,\mu-1}^{\gamma}(\omega(b-a)^\rho)}}\\
&\qquad\times\bigg(1+\frac{\beta(\xi-a)^{\mu-1}E_{\rho,\mu}^{\gamma}(\omega(\xi-a)^\rho)}{(b-a)^{\mu-2}E_{\rho,\mu-1}^{\gamma}(\omega(b-a)^\rho)- \beta (\xi-a)^{\mu-1}E_{\rho,\mu}^{\gamma}(\omega(\xi-a)^\rho)}\bigg),
\end{align*}
we have
\begin{align*}
&x(t)=\frac{(t-a)^{\mu-1}E_{\rho,\mu}^{\gamma}(\omega(t-a)^\rho)}{(b-a)^{\mu-2}E_{\rho,\mu-1}^{\gamma}(\omega(b-a)^\rho)}\\
&\bigg[1+\frac{\beta(\xi-a)^{\mu-1}E_{\rho,\mu}^{\gamma}(\omega(\xi-a)^\rho)}{(b-a)^{\mu-2}E_{\rho,\mu-1}^{\gamma}(\omega(b-a)^\rho)-\beta(\xi-a)^{\mu-1}E_{\rho,\mu}^{\gamma}(\omega(\xi-a)^\rho)}\bigg]\\
&\times \int_{a}^{b}(b-s)^{\mu-2}E_{\rho,\mu-1}^{\gamma}(\omega(b-s)^\rho)q(s)x(s)ds \\
&-\frac{\beta(t-a)^{\mu-1}E_{\rho,\mu}^{\gamma}(\omega(t-a)^\rho)}{(b-a)^{\mu-2}E_{\rho,\mu-1}^{\gamma}(\omega(b-a)^\rho)-\beta(\xi-a)^{\mu-1}E_{\rho,\mu}^{\gamma}(\omega(\xi-a)^\rho)}\\
&\times\int_{a}^{\xi}(\xi-s)^{\mu-1}E_{\rho,\mu}^{\gamma}(\omega(\xi-s)^\rho)q(s)x(s)ds-\int_{a}^{t}(t-s)^{\mu-1}E_{\rho,\mu}^{\gamma}(\omega(t-s)^\rho)q(s)x(s)ds.
\end{align*}
On simplifying,
\begin{align*}
&x(t)=\frac{(t-a)^{\mu-1}E_{\rho,\mu}^{\gamma}(\omega(t-a)^\rho)}{(b-a)^{\mu-2}E_{\rho,\mu-1}^{\gamma}(\omega(b-a)^\rho)}\int_{a}^{t}(b-s)^{\mu-2}E_{\rho,\mu-1}^{\gamma}(\omega(b-s)^\rho)q(s)x(s)ds\\
&+\frac{(t-a)^{\mu-1}E_{\rho,\mu}^{\gamma}(\omega(t-a)^\rho)}{(b-a)^{\mu-2}E_{\rho,\mu-1}^{\gamma}(\omega(b-a)^\rho)}\int_{t}^{b}(b-s)^{\mu-2}E_{\rho,\mu-1}^{\gamma}(\omega(b-s)^\rho)q(s)x(s)ds\\
&+\bigg(\frac{(t-a)^{\mu-1}E_{\rho,\mu}^{\gamma}(\omega(t-a)^\rho)}{(b-a)^{\mu-2}E_{\rho,\mu-1}^{\gamma}(\omega(b-a)^\rho)}\bigg)\bigg[\frac{\beta(t-a)^{\mu-1}E_{\rho,\mu}^{\gamma}(\omega(t-a)^\rho)}{(b-a)^{\mu-2}E_{\rho,\mu-1}^{\gamma}(\omega(b-a)^\rho)-\beta(\xi-a)^{\mu-1}E_{\rho,\mu}^{\gamma}(\omega(\xi-a)^\rho)}\bigg]\\
&\times\int_{a}^{\xi}(b-s)^{\mu-2}E_{\rho,\mu-1}^{\gamma}(\omega(b-s)^\rho)q(s)x(s)ds\\
&+\bigg(\frac{(t-a)^{\mu-1}E_{\rho,\mu}^{\gamma}(\omega(t-a)^\rho)}{(b-a)^{\mu-2}E_{\rho,\mu-1}^{\gamma}(\omega(b-a)^\rho)}\bigg)\\
&\quad\times\bigg[\frac{\beta(t-a)^{\mu-1}E_{\rho,\mu}^{\gamma}(\omega(t-a)^\rho)}{(b-a)^{\mu-2}E_{\rho,\mu-1}^{\gamma}(\omega(b-a)^\rho)-\beta(\xi-a)^{\mu-1}E_{\rho,\mu}^{\gamma}(\omega(\xi-a)^\rho)}\bigg]\\
&\qquad\times\int_{\xi}^{b}(b-s)^{\mu-2}E_{\rho,\mu-1}^{\gamma}(\omega(b-s)^\rho)q(s)x(s)ds\\
&-\frac{\beta(t-a)^{\mu-1}E_{\rho,\mu}^{\gamma}(\omega(t-a)^\rho)}{(b-a)^{\mu-2}E_{\rho,\mu-1}^{\gamma}(\omega(b-a)^\rho)-\beta(\xi-a)^{\mu-1}E_{\rho,\mu}^{\gamma}(\omega(\xi-a)^\rho)}\\
&\qquad\times\int_{a}^{\xi}(\xi-s)^{\mu-1}E_{\rho,\mu}^{\gamma}(\omega(\xi-s)^\rho)q(s)x(s)ds -\int_{a}^{t}(t-s)^{\mu-1}E_{\rho,\mu}^{\gamma}(\omega(t-s)^\rho)q(s)x(s)ds.
\end{align*}
Further, on rearranging the terms, we have
\begin{align*}
x(t)=&\int_{a}^{t}\bigg[\frac{(t-a)^{\mu-1}E_{\rho,\mu}^{\gamma}(\omega(t-a)^\rho)(b-s)^{\mu-2}E_{\rho,\mu-1}^{\gamma}(\omega(b-s)^\rho)}{(b-a)^{\mu-2}E_{\rho,\mu-1}^{\gamma}(\omega(b-a)^\rho)}\\
&-(t-s)^{\mu-1}E_{\rho,\mu}^{\gamma}(\omega(t-s)^\rho)\bigg]q(s)x(s)ds\\
&+\int_{t}^{b}\bigg[\frac{(t-a)^{\mu-1}E_{\rho,\mu}^{\gamma}(\omega(t-a)^\rho)(b-s)^{\mu-2}E_{\rho,\mu-1}^{\gamma}(\omega(b-s)^\rho)}{(b-a)^{\mu-2}E_{\rho,\mu-1}^{\gamma}(\omega(b-a)^\rho)}\bigg]q(s)x(s)ds\\
&+\frac{\beta(t-a)^{\mu-1}E_{\rho,\mu}^{\gamma}(\omega(t-a)^\rho)}{(b-a)^{\mu-2}E_{\rho,\mu-1}^{\gamma}(\omega(b-a)^\rho)-\beta(\xi-a)^{\mu-1}E_{\rho,\mu}^{\gamma}(\omega(\xi-a)^\rho)}\\
&\quad\times\int_{a}^{\xi}\bigg[\frac{(\xi-a)^{\mu-1}E_{\rho,\mu}^{\gamma}(\omega(\xi-a)^\rho)(b-s)^{\mu-2}E_{\rho,\mu-1}^{\gamma}(\omega(b-s)^\rho)}{(b-a)^{\mu-2}E_{\rho,\mu-1}^{\gamma}(\omega(b-a)^\rho)}\\
&-(\xi-s)^{\mu-1}E_{\rho,\mu}^{\gamma}(\omega(\xi-s)^\rho)\bigg]q(s)x(s)ds\\
&+\bigg(\frac{(t-a)^{\mu-1}E_{\rho,\mu}^{\gamma}(\omega(t-a)^\rho)}{(b-a)^{\mu-2}E_{\rho,\mu-1}^{\gamma}(\omega(b-a)^\rho)}\bigg)\\
&\quad\times\bigg[\frac{\beta(\xi-a)^{\mu-1}E_{\rho,\mu}^{\gamma}(\omega(\xi-a)^\rho)}{(b-a)^{\mu-2}E_{\rho,\mu-1}^{\gamma}(\omega(b-a)^\rho)-\beta(\xi-a)^{\mu-1}E_{\rho,\mu}^{\gamma}(\omega(\xi-a)^\rho)}\bigg]\\
&\qquad\times \int_{\xi}^{b}(b-s)^{\mu-2}E_{\rho,\mu-1}^{\gamma}(\omega(b-s)^\rho)q(s)x(s)ds,
\end{align*}
therefore,
\begin{align*}
&x(t)=\int_{a}^{b}G(t,s)q(s)x(s)ds\\
&+\frac{\beta(t-a)^{\mu-1}E_{\rho,\mu}^{\gamma}(\omega(t-a)^\rho)}{(b-a)^{\mu-2}E_{\rho,\mu-1}^{\gamma}(\omega(b-a)^\rho)-\beta(\xi-a)^{\mu-1}E_{\rho,\mu}^{\gamma}(\omega(\xi-a)^\rho)}
\int_{a}^{b}G(\xi,s)q(s)x(s)ds,
\end{align*}
where the Green's function $G(t,s)$ is as in \eqref{b4}.
\begin{thm}
The Green's function \eqref{b4} satisfies the following properties:\\
$(a)$ $G(t,s)\geq0,$ for all $(t,s)\in[a,b]\times[a,b]$;\\
$(b)$ $G(t,s)$ is nondecreasing function with respect to the first variable;\\
$(c)$ $0\leq G(a,s)\leq G(t,s)\leq G(b,s), \qquad(t,s)\in [a,b]\times[a,b].$
\end{thm}
\textbf{proof (a)}.
We set two function as
\begin{align*}
g_{1}(t,u)=&\frac{(t-a)^{\mu-1}E_{\rho,\mu}^{\gamma}(\omega(t-a)^\rho){(b-s)^{\mu-2}E_{\rho,\mu-1}^{\gamma}(\omega(b-s)^\rho)}}{{(b-a)^{\mu-2}E_{\rho,\mu-1}^{\gamma}(\omega(b-a)^\rho)}}\\
&-(t-s)^{\mu-1}E_{\rho,\mu}^{\gamma}(\omega(t-s)^\rho),\quad \quad a\leq s\leq t\leq b,\\
\end{align*}
and \\
\begin{align*}
g_{2}(t,u)=&\frac{(t-a)^{\mu-1}E_{\rho,\mu}^{\gamma}(\omega(t-a)^\rho){(b-s)^{\mu-2}E_{\rho,\mu-1}^{\gamma}(\omega(b-s)^\rho)}}{{(b-a)^{\mu-2}E_{\rho,\mu-1}^{\gamma}(\omega(b-a)^\rho)}},\quad  \quad a\leq t\leq s\leq b.
\end{align*}
It is clear that $g_{2}(t,u)\geq0$. So to prove $(a)$, we should show that $g_{1}(t,u)\geq0$, or equivalently\\
\begin{align*}
&\frac{(t-a)^{\mu-1}E_{\rho,\mu}^{\gamma}(\omega(t-a)^\rho){(b-s)^{\mu-2}E_{\rho,\mu-1}^{\gamma}(\omega(b-s)^\rho)}}{{(b-a)^{\mu-2}E_{\rho,\mu-1}^{\gamma}(\omega(b-a)^\rho)}}
\geq(t-s)^{\mu-1}E_{\rho,\mu}^{\gamma}(\omega(t-s)^\rho),\\
&\text{therefore it is sufficient to prove that}\\
&(i)\frac{(t-a)^{\mu-1}(b-s)^{\mu-2}}{(b-a)^{\mu-1}}\geq(t-s)^{\mu-1}\\
&(ii)\frac{E_{\rho,\mu}^{\gamma}(\omega(t-a)^\rho)E_{\rho,\mu-1}^{\gamma}(\omega(b-s)^\rho)}{E_{\rho,\mu-1}^{\gamma}(\omega(b-a)^\rho)}\geq E_{\rho,\mu}^{\gamma}(\omega(t-s)^\rho)\\
\end{align*}
 Consider,
\begin{align*}
 &\frac{(t-a)^{\mu-1}(b-s)^{\mu-2}}{(b-a)^{\mu-1}}-(t-s)^{\mu-1}\\
 &=(t-a)^{\mu-1}\bigg(\frac{b-s}{b-a}\bigg)^{\mu-2}-((t-a)-(s-a))^{\mu-1}\\
 &= (t-a)^{\mu-1}\bigg(\frac{(b-a)-(s-a)}{b-a}\bigg)^{\mu-2}-(t-a)^{\mu-1}\bigg(1-\frac{s-a}{t-a}\bigg)^{\mu-1}\\
 &= (t-a)^{\mu-1}(1-\frac{s-a}{b-a})^{\mu-2}-(t-a)^{\mu-1}\bigg(1-\frac{s-a}{t-a}\bigg)L^{\mu-1}\\
 &\geq(t-a)^{\mu-1}(1-\frac{s-a}{b-a})^{\mu-2}-(t-a)^{\mu-1}\bigg(1-\frac{s-a}{b-a}\bigg)^{\mu-1}\quad(\because t\leq b)\\
 &=(t-a)^{\mu-1}\bigg(1-\frac{s-a}{b-a}\bigg)^{\mu-2}\bigg(\frac{s-a}{b-a}\bigg)\geq0.
\end{align*}
Hence (i) is proved. For the proof of (ii) refer (Theorem $2$, in \cite{EA})\\
\textbf{proof (b).} Proof of this is similar to (Theorem $2$, in \cite{EA}) for $\mu=\mu-1$.\\
\textbf{Proof (c).} Proof of this follows from \textbf(b).
\begin{thm}
Suppose that problem \eqref{3} has a nontrivial continuous solution, then
\begin{align*}
&\int_{a}^{b} \bigg[\frac{(b-a)^{\mu-1}E_{\rho,\mu}^{\gamma}(\omega(b-a)^\rho)(b-s)^{\mu-2}E_{\rho,\mu-1}^{\gamma}(\omega(b-s)^\rho)}{(b-a)^{\mu-2}E_{\rho,\mu-1}^{\gamma}(\omega(b-a)^\rho)}\\
&\qquad-(b-s)^{\mu-1}E_{\rho,\mu}^{\gamma}(\omega(b-s)^\rho)\bigg]|q(s)|ds\\
&\geq \bigg(1+\frac{\beta(\xi-a)^{\mu-1}E_{\rho,\mu}^{\gamma}(\omega(\xi-a)^\rho)}{(b-a)^{\mu-2}E_{\rho,\mu-1}^{\gamma}(\omega(b-a)^\rho)-\beta (\xi-a)^{\mu-1}E_{\rho,\mu}^{\gamma}(\omega(\xi-a)^\rho)}\bigg)^{-1}
\end{align*}
\end{thm}
\textbf{Proof}.
Consider the Banach space\\
$C[a,b]=\{u:[a,b]\longrightarrow\mathbb{R}\setminus \text{$u$ is continuous}\}$\\
$\|u\|_{\infty}= max\{|u(t)|: a\leq t \leq b\},u\in C[a,b]$.\\
By theorem \ref{T1}, a solution $x\in C[a,b]$ of \eqref{3} has the expression
\begin{align*}
&x(t)=\int_{a}^{b}G(t,s)q(s)x(s)ds\\
&+\frac{\beta(t-a)^{\mu-1}E_{\rho,\mu}^{\gamma}(\omega(t-a)^\rho)}{(b-a)^{\mu-2}E_{\rho,\mu-1}^{\gamma}(\omega(b-a)^\rho)-\beta (\xi-a)^{\mu-1}E_{\rho,\mu}^{\gamma}(\omega(\xi-a)^\rho)}
&\int_{a}^{b}G(\xi,s)q(s)x(s)ds,\quad a\leq t\leq b.
\end{align*}
From this, for any $t\in[a,b]$, we have
\begin{align*}
|x(t)|\leq  &\|x\|_{\infty}\int_{a}^{b}|G(t,s)||q(s)|ds \\
&+\frac{\beta(t-a)^{\mu-1}E_{\rho,\mu}^{\gamma}(\omega(t-a)^\rho)\|x\|_{\infty}}{(b-a)^{\mu-2}E_{\rho,\mu-1}^{\gamma}(\omega(b-a)^\rho)-\beta (\xi-a)^{\mu-1}E_{\rho,\mu}^{\gamma}(\omega(\xi-a)^\rho)}
&\int_{a}^{b}|G(\xi,s)||q(s)|ds,
\end{align*}
therefore,
\begin{align*}
|x(t)|\leq &\|x\|_{\infty}\int_{a}^{b}|G(b,s)||q(s)|ds\\
&+\frac{\beta(b-a)^{\mu-1}E_{\rho,\mu}^{\gamma}(\omega(b-a)^\rho)\|x\|_{\infty}}{(b-a)^{\mu-2}E_{\rho,\mu-1}^{\gamma}(\omega(b-a)^\rho)-\beta (\xi-a)^{\mu-1}E_{\rho,\mu}^{\gamma}(\omega(\xi-a)^\rho)}
&\int_{a}^{b}|G(b,s)||q(s)|ds,
\end{align*}
which yields
\begin{align*}
\|x\|_{\infty}&\leq\|x\|_{\infty}\bigg(1+\frac{\beta(b-a)^{\mu-1}E_{\rho,\mu}^{\gamma}(\omega(b-a)^\rho)}{(b-a)^{\mu-2}E_{\rho,\mu-1}^{\gamma}(\omega(b-a)^\rho)-\beta(\xi-a)^{\mu-1}E_{\rho,\mu}^{\gamma}(\omega(\xi-a)^\rho)}\bigg)
\int_{a}^{b}G(b,s)|q(s)|ds,\\
\end{align*}
As $x$ is a nontrivial, we have
\begin{align*}
&1\leq \bigg(1+\frac{\beta(b-a)^{\mu-1}E_{\rho,\mu}^{\gamma}(\omega(b-a)^\rho)}{(b-a)^{\mu-2}E_{\rho,\mu-1}^{\gamma}(\omega(b-a)^\rho)-\beta (\xi-a)^{\mu-1}E_{\rho,\mu}^{\gamma}(\omega(\xi-a)^\rho)}\bigg)
\int_{a}^{b}|G(b,s)||q(s)|ds.\\
\end{align*}
\begin{align*}
\int_{a}^{b}|G(b,s)||q(s)|ds\geq\bigg(1+\frac{\beta(b-a)^{\mu-1}E_{\rho,\mu}^{\gamma}(\omega(b-a)^\rho)}{(b-a)^{\mu-2}E_{\rho,\mu-1}^{\gamma}(\omega(b-a)^\rho)-\beta (\xi-a)^{\mu-1}E_{\rho,\mu}^{\gamma}(\omega(\xi-a)^\rho)}\bigg)^{-1},
\end{align*}
therefore
\begin{align*}
\int_{a}^{b} \bigg[\frac{(b-a)^{\mu-1}E_{\rho,\mu}^{\gamma}(\omega(b-a)^\rho)(b-s)^{\mu-2}E_{\rho,\mu-1}^{\gamma}(\omega(b-s)^\rho)}{(b-a)^{\mu-2}E_{\rho,\mu-1}^{\gamma}(\omega(b-a)^\rho)}
-(b-s)^{\mu-1}E_{\rho,\mu}^{\gamma}(\omega(b-s)^\rho)\bigg]|q(s)|ds\\
\geq \bigg(1+\frac{\beta(\xi-a)^{\mu-1}E_{\rho,\mu}^{\gamma}(\omega(\xi-a)^\rho)}{(b-a)^{\mu-2}E_{\rho,\mu-1}^{\gamma}(\omega(b-a)^\rho)-\beta (\xi-a)^{\mu-1}E_{\rho,\mu}^{\gamma}(\omega(\xi-a)^\rho)}\bigg)^{-1} \\
\end{align*}
Hence the result.


\begin{thebibliography}{9}
\bibitem {ND} N. Abuj and D. Pachpatte, \emph{Lyapunov type inequality for fractional differential equation with $k$-Prabhakar derivative},
Classical Analysis and ODEs (math.CA) (2017), arXiv:1702.01562v1 [math.CA], 1-11.
\bibitem{BH} R. Brown, D. Hinton, \emph{Opial's inequality and oscillation of $2$nd order equations}, \textit{Proc. Amer. Math. Soc.} 125 (1997), 1123-1129.
\bibitem{CS1} S. Cheng, \emph{A discrete analouge of the inequality of Lyapunov}, \textit{Hokkaido Math. J.} 12 (1983), 105-112.
\bibitem{CS2} S. Cheng, \emph{Lyapunov type inequalities for differential and difference equations}, \textit{Fasc. Math.} 23 (1991), 25-41.
\bibitem{CD} D. Cakmak, \emph{Lyapunov-type integral inequalities for certain higher order differential equations}, \textit{Appl. Math. Comput.} 216 (2010), 368-373.
\bibitem{CSS} I. Cabrera, K. Sadarangani and B. Samet, \emph{Hartman-Wintner-type inequalities for a class of nonlocal fractional boundary value problems}, \textit{Math. Methods Appl. Sci.} 40 (2017), 129-136.
\bibitem{EA} S. Eshaghi and A. Ansari, \emph{Lyapunov inequality for fractional differential equations with Prabhakar derivative}, \textit{Math. Inequal. Appl.} 19 (2016), 349-358.
\bibitem{FR1} R. Ferreira, \emph{A Lyapunov-type inequality for a fractional boundary value problem}, \textit{Fract. Calc. Appl. Anal.} 16(2013), 978-984.
\bibitem{RAC2} R. Ferriera, \emph{On a Lyapunov-type inequality and the zeros of a certain Mittag-Leffler function}, \textit{J. Math. Anal. Appl.} 412 (2014), 1058-1063.
\bibitem{RG} R. Garra, R. Gorenflo, F. Polito and Z. Tomovski, \emph{Hilfer Prabhakar derivative and some applications}, \textit{Appl. Math. Comput.} 242 (2014), 576-589.
\bibitem{HW} P. Hartman, A. Wintner, \emph{On an oscillation criterion of Lyapunov}, \textit{Amer. J. Math.} 73 (1951), 885-890.
\bibitem{MJ1} M. Jleli and B. Samet, \emph{Lyapunov-type inequalities for a fractional differential equation with mixed boundary conditions}, \textit{Math. Inequal. Appl.} 18 (2015), 443-451.
\bibitem{MJ2} M. Jleli and B. Samet, \emph{Lyapunov-type inequalities for fractional boundary value problems}, \textit{Electr. j. differ. equ.} 88 (2015), 1-11.
\bibitem{AML} A. Lyapunov, \emph{Probl$\acute{e}$me g$\acute{e}$n$\acute{e}$ral de la stabilit$\acute{e}$ du mouvement}, \textit{Ann. Fac. Sci. Toulouse Math.} 2 (1907), 203-407.
\bibitem{LYH} C. Lee, C. Yeh, C. Hong, R. Agarwal \emph{Lyapunov and Wirtinger inequalities}, \textit{Appl. Math. Lett.} 17 (2004), 847-853.
\bibitem{DNA} D. Pachpatte, N. Abuj and A. Khandagale, \emph{Lyapunov type inequality for hybrid fractional differential equation with Prabhakar derivative}, \textit{Int. J. Pure Appl. Math.} 113 (2017), 563-574.
\bibitem{PB} B. Pachpatte, \emph{On Lyapunov-type inequalities for certain higher order differential equations}, \textit{J. Math. Anal. Appl.} 195 (1995), 527-536.
\bibitem{TRP} T. R. Prabhakar, \emph{A Singular integral equation with a generalized Mittag-Leffler function in the kernel}, \textit{Yokohama Math. J.} 19(1971), 7-15.
\bibitem{PP} N. Parhi, S. Panigrahi, \emph{On Lyapunov-type inequality for third order differential equations}, \textit{J. Math. Anal. Appl.} 233 (1999), 445-460.
\bibitem{RB}  J. Rong, C. Bai, \emph{Lyapunov-type inequality for a fractional differential equation with fractional boundary condition}, \textit{Adv. Difference Equ.} 82 (2015), 1-10.
\bibitem{RS} O. Regan, B. Samet, \emph{Lyapunov-type inequality for a class of fractional differential equations}, \textit{J. Inequal. Appl.} 247 (2015), 1-10.
\bibitem {TUC} A. Tiryaki, M. Unal, D. Cakmak, \emph{Lyapunov-type inequalities for nonlinear systems}, \textit{Math. Anal. Appl.} 332 (2007) 497-511.
\bibitem{YLL} X. Yang, K. Lo, \emph{Lyapunov-type inequality for a class of even-order differential equations}, \textit{Appl. Math. Comput.} 215 (2010), 3884-3890.
\end{thebibliography}
\end{document}